\documentclass[10pt]{amsart}
\usepackage[cp1251]{inputenc}
\usepackage[english,russian]{babel}
\usepackage{amsmath}
\usepackage{amssymb}
\usepackage{amsfonts}
\usepackage{graphicx}

\newtheorem{theorem}{Theorem}

\begin{document}
\sloppy
\renewcommand{\refname}{References}
\renewcommand{\proofname}{Proof.}
\renewcommand{\figurename}{Fig.}

\thispagestyle{empty}

\title{Asymptotics of joint orderings of compound Poisson fields}
\author{{Mikhail Chebunin, Artyom Kovalevskii}}
\address{Mikhail Georgievich Chebunin
\newline\hphantom{iii} Karlsruhe Institute of Technology, 
\newline\hphantom{iii} Institute of Stochastics, 
\newline\hphantom{iii} 76131, Karlsruhe, Germany;
\newline\hphantom{iii} Novosibirsk State University, 
\newline\hphantom{iii} Pirogova str., 1,
\newline\hphantom{iii} 630090, Novosibirsk, Russia}

\email{chebuninmikhail@gmail.com}

\address{Artyom Pavlovich Kovalevskii
\newline\hphantom{iii} Novosibirsk State Technical University, 
\newline\hphantom{iii} pr. K. Marksa, 20,
\newline\hphantom{iii} 630073, Novosibirsk, Russia;
\newline\hphantom{iii} Novosibirsk State University, 
\newline\hphantom{iii} Pirogova str., 1,
\newline\hphantom{iii} 630090, Novosibirsk, Russia}

\email{artyom.kovalevskii@gmail.com}

\thanks{\sc Chebunin, M.G., Kovalevskii, A.P.,
Asymptotics of joint orderings of compound Poisson fields}
\thanks{\copyright \ 2022 Chebunin M.G., Kovalevskii A.P}
\thanks{\rm The work is supported by Mathematical Center in Akademgorodok 
under agreement No. 075-15-2019-1675 with the Ministry of Science and Higher Education of the Russian Federation.
}

\maketitle {\small
\begin{quote}
\noindent{\sc Abstract. }  
We are developing a new method for the analysis of queuing systems with heterogeneous in time and space
compound (marked) Poisson input flow. The state space of the input flow is embedded in a higher-dimensional space
with a homogeneous marked Poisson field on it. We prove
limit theorems for partial sums of marks under the ordering of field points by coordinates.

\medskip

 \end{quote}
}

\section{Introduction}

Modern queuing systems, including cloud services, deal with a time-varying flow of customers, see Bouterse and Perros (2012). New methods of analysis and optimization take into account not only queues on servers, but also the location of users on the surface of the earth and even in space. Each client has a number of characteristics that can be represented by a vector.
The heterogeneity of the arrival of clients in time and space can be leveled by the introduction of additional
dimensions. So, if the intensity $\lambda(t)$ changes with time $ t \in [0,T]$, then such a marked (compound) inhomogeneous Poisson flow corresponds to a two-dimensional homogeneous Poisson field on the set
$\{0 \le x_1 \le T, \ 0 \le x_2 \le \lambda(x_1)\}$, see Section 3.

We prove Theorem 1 on convergence of a random field of the sums of marks ${\bf Y}_i$ of those points ${\bf X}_i$ of the Poisson process, which satisfy the coordinate-wise inequality ${\bf X}_i \le {\bf u}$ for ${\bf u} \in \mathbb{R}^{d_1}$.

Then we prove Theorem 2 on convergence of processes of partial sums of vectors composed of marks that
are ordered in accordance with the increase of each of the coordinates of the points of the Poisson field.

Theorem 3 describes the limiting behavior of a particular model of a compound Poisson process that is inhomogeneous in time when its points are ordered in time and intensity.

The method of proof is based on a Poissonization of the results of Chebunin and Kovalevskii (2021),
which, in turn, are based on a partial generalization to the multidimensional case of the results of
Davydov and Egorov (2000). Both of these papers use a very general theorem of Ossiander (1987).

The rest of the paper is organized as follows. Section 2 contains general results for orderings of 
 compound Poisson fields. Section 3 describes a concrete model of marked non-homogeneous Poisson process
 and the limiting result for it. Proofs are presented in Section 4.

\section{Results for compound Poisson fields}

Let $A \subset \mathbb{R}^{d_1}$ be a compact Borel set with the Lebesgue measure
$0<\mu(A)<\infty$, and there is a compound Poisson field in $A$ with intensity $\nu>0$. That is, the number 
$\eta=\eta_{\nu}$ of points ${\bf X}_i$ in $A$ has
Poisson distribution with  parameter $\nu\mu(A)$, the points are distributed uniformly in $A$
under the condition of fixed $\eta$.

Any point ${\bf X}_i$ is associated with a mark ${\bf Y}_i$ that takes values in ${\bf R}^{d_2}$, $0<d_2<\infty$.
Vectors $({\bf X}_i, {\bf Y}_i)$ are mutually independent and identically distributed copies of a random vector 
$({\bf X}, {\bf Y})$ such that 
${\bf X}=(X^{(1)},\ldots, X^{(d_1)} )$ takes values in $A$, ${\bf Y}$ takes values in $\mathbb{R}^{d_2}$. 

Matrix $(X,Y)$ has a random number $\eta$ of rows $({\bf X}_i, {\bf Y}_i)$.

Denote $X^{(k)}_{\eta, 1} \leq X^{(k)}_{\eta, 2} \leq \cdots \leq
X^{(k)}_{\eta, \eta}$, $1\le k\le d_1$, the order statistics of the $k$-th column of matrix $X$, and 
${\bf Y}^{(k)}_{\eta, 1}, {\bf Y}^{(k)}_{\eta, 2}, \ldots, {\bf Y}^{(k)}_{\eta, \eta}$ 
the corresponding values of the vectors 
${\bf Y}_i$. The random vectors $\left({\bf Y}^{(k)}_{\eta, i}, i \leq \eta \right)$ are called induced order statistics (concomitants).

We look into the asymptotic behavior of the random field
$$
\begin{aligned}
{\bf Q}_{\nu}({\bf u}) &= \sum_{j=1}^{\eta_{\nu}} {\bf Y}_{ j} \mathbf{1}\left({\bf X}_{ j}\le {\bf u}\right)
= \sum_{j=1}^{\eta_{\nu}} {\bf Y}_{ j} \mathbf{1}\left({X}_{ j}^{(1)}\le u^{(1)}, \ldots, {X}_{ j}^{(d_1)}\le u^{(d_1)}\right), \ {\bf u}\in {\bf R}^{d_1}
\end{aligned}
$$
as $\nu\to \infty$.

Using the asymptotics of ${\bf Q}_{\nu}({\bf u})$, we study the asymptotics of $d_1\times d_2$-dimensional process 
of sums of induced order statistics under different orderings
$$
\begin{aligned}
{\bf Z}_{\nu}(t) &= \left(
\sum_{j=1}^{[\eta_{\nu} t]} {\bf Y}^{(1)}_{\eta_{\nu}, j},
\sum_{j=1}^{[\eta_{\nu} t]} {\bf Y}^{(2)}_{\eta_{\nu}, j},
 \dots,
 \sum_{j=1}^{[\eta_{\nu} t]} {\bf Y}^{(d_1)}_{\eta_{\nu}, j} 
 \right), \ t \in [0, 1].
\end{aligned}
$$

Let ${\bf m}({\bf v})={\bf E}({\bf Y} \mid {\bf X}={\bf v})$, ${\bf v} \in  A$, and ${\bf f}({\bf u})=\int_{{\bf -\infty}}^{\bf u} {\bf  m}({\bf v}) 
{\bf 1}({\bf v}\in A) d {\bf v}$.

 Let
$$
\sigma^{2}({\bf v})={\bf E}\left\{({\bf Y}-{\bf m}({\bf X}))^{T}({\bf Y}-{\bf m}({\bf X})) \mid {\bf X}={\bf v}\right\}
$$
be the conditional covariance matrix of ${\bf Y}$ and $\sigma({\bf v})$ be the positive definite matrix such that $\sigma({\bf v})^{T} \sigma({\bf v})=\sigma^{2}({\bf v})$.

All our limit fields and processes are continuous, so we use the uniform metric. 
Let $\|\cdot\|$ denote the Euclidean norm in the corresponding space.
Our random field ${\bf Q}_{\nu}$ takes values in the space $B(A; \, {\bf R}^{d_2})$
of bounded measurable functions with the Borel $\sigma$-algebra $\mathcal{B}$. This space is not separable, so we note that the random field 
${\bf Q}_{\nu}$
takes values in its subset $D$ with the smaller $\sigma$-algebra $\mathcal{D}$. This $\sigma$-algebra 
is generated by the $d_1$-dimensional analog of  Skorohod metrics, see Straf (1972).
So, let $D$ be the uniform closure, in the space $B(A; \, {\bf R}^{d_2})$, of the vector subspace
of simple functions (that is, linear combinations of step functions).

For ${\bf x},{\bf y} \in D$, let the ``Skorohod'' distance be
\[
d({\bf x},{\bf y})=\inf\{\min(||{\bf x}-\mathbf{y\lambda}||_S,||\lambda||_S): \ \lambda \in \Lambda\},
\]
where 
\[
||{\bf x}-{\bf y}\lambda||_S=\sup\{|{\bf |x}({\bf t})-{\bf y}(\lambda({\bf t}))||, \ {\bf t} \in A\},
\] 
$||\lambda||_S=\sup\{||\lambda({\bf t})-{\bf t}||, \ {\bf t} \in A\}$, $\Lambda$ is 
the group of all transformations $\lambda$ of the form $\lambda(t_1,\ldots,t_{d_1})=
(\lambda_1(t_1), \ldots, \lambda_{d_1}(t_{d_1}))$, where each $\lambda_i:[0,1] \to [0,1]$ is continuous, 
strictly increasing, and fixes 0 and 1. 

With respect to the corresponding metric topology, $D$ is complete and separable, and its Borel $\sigma$-algebra
$\mathcal{D}$ coincides with the $\sigma$-algebra generated by coordinate mappings, see Section 3 of Bickel $\&$
 Wichura (1971) for details. 
 
 The process ${\bf Q}_{\nu}$ takes values in ${D}$. So, it is $\mathcal{D}$-measurable.

We define the weak convergence as follows (cf.
Dudley, 1967): $\lim_{\nu\to\infty} {\bf E} f({\bf Q}_{\nu}) = {\bf E} f({\bf Q})$ for every bounded, continuous, and
$\mathcal{D}$-measurable function $f:\, D \to {\bf R}^{d_2}$.

Davydov $\&$ Zitikis (2008) in their Proposition 1 have proved that the limitation of the class of functions $f$ to the uniform continuity (instead of the continuity) gives the same definition of the weak convergence (that is, these definitions are equivalent).

We use the symbol $\Rightarrow$ to denote the weak convergence of random fields in the sense that has 
been mentioned above. We use the same symbol for the weak convergence of random variables and the weak convergence of stochastic processes in the uniform topology.

The following Theorem 1 propagates the result of Lemma 1 by Chebunin and Kovalevskii (2021) to compound Poisson fields.

\begin{theorem}
If ${\bf E} \|{\bf Y}\|^2 < \infty$ then
$\widetilde{\bf Q}_{\nu}=\frac{{\bf Q}_{\nu} - \eta_{\nu}{\bf f}}{\sqrt{\eta_{\nu}}} \Rightarrow {\bf Q}$, a centered Gaussian field with covariance 
\[
K({\bf u}_1, {\bf u}_2) = {\bf E} {\bf Q}^T({\bf u}_1) {\bf Q}({\bf u}_2) = 
\int\limits_{A_{{\bf u}_1,{\bf u}_2}} \sigma^2({\bf v})  d {\bf v} 
\]
\[
+ 
\int\limits_{A_{{\bf u}_1,{\bf u}_2}} {\bf m}^T({\bf v}) {\bf m}({\bf v})   d {\bf v}
- 
\int\limits_{A_{{\bf u}_1,{\bf u}_2}} {\bf m}^T({\bf v})  d {\bf v}
\int\limits_{A_{{\bf u}_1,{\bf u}_2}} {\bf m}({\bf v})   d {\bf v},
\]
${A_{{\bf u}_1,{\bf u}_2}}=\{{\bf v} \in {\bf R}^{d_1}: \ {\bf v } \in A, {\bf v} \le {\bf u}_1, {\bf v} \le{\bf u}_2 \}$.
\end{theorem}

The following Theorem 2 propagates the result of Lemma 2 by Chebunin and Kovalevskii (2021) to compound Poisson fields.

\begin{theorem}
 If ${\bf E} \|{\bf Y}\|^2 < \infty$, ${\bf m}\equiv {\bf 0}$ then $\widetilde{\bf Z}_{\nu}
 =\frac{{\bf Z}_{\nu}}{\sqrt{\eta_{\nu}}} \Rightarrow {\bf Z}$, a centered Gaussian $(d_1\times d_2)$-dimensional process with
covariance matrix function  ${\bf E}{\bf Z}^T(t_1){\bf Z}(t_2)=\widetilde{K}(t_1,t_2) =(\widetilde{K}_{ij} (t_1,t_2)_{i,j=1}^{d_1}$, 
\[
\widetilde{K}_{ij} (t_1,t_2) = 
\int\limits_{B_{i,t_1,j,t_2}}  \sigma^2({\bf v})  d {\bf v},
\]
$B_{i,t_1,j,t_2}=\{{\bf v} \in {\bf R}^{d_1}: \ {\bf v} \in A, v_i \le t_1, v_j \le t_2 \}$.
\end{theorem}

\section{Multiple ordering of non-homogeneous Poisson process}

We interpret the results of the previous section for non-homogeneous Poisson processes.
So we take set $A$ of special form,
\[
A=\{{\bf x}=(x_1,x_2): \ 0 \le x_1 \le T, 0 \le x_2 \le \lambda(x_1)\},
\]
$T>0$, $\lambda=\{\lambda(t), 0 \le t \le T\}>0$ is the positive Borel function on $[0,T]$. 

We interpret $x_1$ as time, $\lambda(x_1)$ as the Poisson parameter at time $x_1$.

So $({\bf X}_i,{\bf Y}_i)=((X_{i1},X_{i2}),{\bf Y}_i)$.

We suppose that ${\bf m}(x_1,x_2)$ does not depend on $x_2$:
Let
\[
{\bf m}(x_1,x_2)={\bf m}(x_1)={\bf E}({\bf Y}_1|X_{11}=x_1).
\]

We order vectors $({\bf X}_1,{\bf Y}_1), \ldots, ({\bf X}_{\eta},{\bf Y}_{\eta}) $, $\eta=\eta_{\nu}$,  using two different algorithms.

The first one: $({\bf X}_i^*,{\bf Y}_i^*)=(X_{i1}^*,X_{i2}^*,{\bf Y}_i^*)$
are ordered by time, that is,
\[
X_{11}^*< X_{21}^*<\ldots<X_{\eta 1}^*.
\]

The second one: $({\bf X}_i^@,{\bf Y}_i^@)=(X_{i1}^@,X_{i2}^@,{\bf Y}_i^@)$
are ordered by the corresponding intensity, that is,
\[
\lambda(X_{11}^@)\le \lambda(X_{21}^@)\le\ldots\le\lambda(X_{\eta 1}^@).
\]

If $\lambda(X_{i1}^@)= \lambda(X_{j1}^@)$ then the order of corresponding pairs is random.

So if $\lambda\equiv const$ then all the pairs $({\bf X}_i^@,{\bf Y}_i^@)$ are in a random order,
this case corresponds to the homogeneous Poisson process.

\begin{theorem} If $0<\sigma^2({\bf x})<\infty$ for any ${\bf x} \in A$ then
\[
\left\{
\eta_{\nu}^{-1/2}\sum_{i=1}^{[\eta_{\nu}t]} \left( {\bf Y}_i^* -{\bf  m}(X_{i1}^*),
 {\bf Y}_i^@ - {\bf m}(X_{i1}^@)\right),
\ 0 \le t \le 1,
\right\}
\]
converges weakly in uniform metrics as $\nu \to \infty$ to a centered Gaussian $2d_2$-dimensional process 
$\left\{({\bf V}_1(t),{\bf V}_2(t)),0\le t \le 1\right\}$ with covariances
\[
{\bf E} {\bf V}_1^T(t_1){\bf V}_2(t_2)=
 \int_{G_{12}(t_1,t_2)} \sigma^2({\bf x})  \, d {\bf x},
\]
$G_{12}(t_1,t_2)=\{(x_1,x_2): \ 0\le x_1\le t_1, \lambda(x_1)\le \lambda(t_2), 0 \le x_2 \le \lambda(x_1)\}$,

\[
{\bf E} {\bf V}_1^T(t_1){\bf V}_1(t_2)=
 \int_{G_{11}(t_1,t_2)} \sigma^2({\bf x})  \, d {\bf x},
\]
$G_{11}(t_1,t_2)=\{(x_1,x_2): \ 0\le x_1\le \min (t_1, t_2), 0 \le x_2 \le \lambda(x_1)\}$,

\[
{\bf E} {\bf V}_2^T(t_1){\bf V}_2(t_2)=
 \int_{G_{22}(t_1,t_2)} \sigma^2({\bf x})  \, d {\bf x},
\]
$G_{22}(t_1,t_2)=\{(x_1,x_2): \lambda(x_1)\le \min (\lambda(t_1),\lambda(t_2)) , 0 \le x_2 \le \lambda(x_1)\}$.
\end{theorem}

\section{Proofs}

{\em Proof of Theorem 1}

Note that a bounded set $A$ can be placed in a $d_1$-dimensional cube. Next, we consider a compound Poisson field with intensity measure $\nu\mu$ over this cube, where $\nu$ is a positive number and $\mu()$ is the Lebesgue measure. As before, we will associate each point of the Poisson field that falls into the set $A$ with a random mark ${\bf Y}_i$, and the points that do not fall into the set $A$ with zero. Thus, using the splitting property of the Poisson field into independent fields, we can assume without loss of generality that the set $A$ is a $d_1$-dimensional cube with side $a>0$ and vertex ${\bf x}_0$ (the smallest point). It is clear that any point ${\bf x} \in A$ of this cube can be easily translated (one-to-one) to the point ${\bf x}'=({\bf x}-{\bf x}_0)/a$ from $d_1$-dimensional cube with side $1$ and vertex at zero.

Note that for a fixed $\eta=N$, the random field ${\bf Q}'_{N}({\bf x}')$ satisfies conditions of Lemma 1 in Chebunin and Kovalevskii (2021), where
\[
{\bf Q}'_N({\bf x}')= \frac{\sum_{i=1}^N {\bf Y}_i {\bf 1}(({\bf X}_i-{\bf x}_0)/a \le {\bf x}') - N \, {\bf E} 
\left({\bf Y}_1 {\bf 1}(({\bf X}_1-{\bf x}_0)/a \le {\bf x}')\right)}{\sqrt{N}}.
\]

Using this lemma for any measurable continuous functional $g$ we have
\[
{\bf E} (g({\bf Q}'_{\eta_{\nu}}({\bf x}'))-g({\bf Q}'({\bf x}'))
={\bf E} ({\bf E} (g({\bf Q}'_{\eta_{\nu}}({\bf x}'))-g({\bf Q}'({\bf x}'))|{\eta_{\nu}}))
\stackrel{def}{=}{\bf E}(q_{{\eta_{\nu}}})  \to 0
\]
as $\nu \to \infty$. Really, from the Lemma we have  $q_n \to 0$ as $n\to\infty$, so sequence $q_n$ 
is a fundamental one. And  $\eta_{\nu}\to \infty$ a.s. as $\nu\to\infty$, so 
\[
{\bf E}(q_{\eta_{\nu}})=\sum_{k=0}^\infty q_k {\bf P}(\eta_{\nu}=k)\pm q_{k_0}
\le {\bf P}(\eta_{\nu}\le k_0) \max_{k\le k_0}(|q_k-q_{k_0}|) +q_{k_0}+\sup_{k>k_0} |q_k-q_{k_0}|\to 0
\]
as $\nu\to\infty$ and $k_0\to\infty$.

The proof is complete.

{\em Proof of Theorem 2}

Theorem 1 is the Poisson analog of Lemma 1 in Chebunin and Kovalevskii (2021). So the proof of Theorem 2 goes by lines 
of Lemma 2  in Chebunin and Kovalevskii (2021) using SLLN for 
$\eta_{\nu}$. We omit the details.

{\em Proof of Theorem 3} 

Note that our assumptions entail zero expectations of  ${\bf Y}_i^* -{\bf  m}(X_{i1}^*)$ and
 ${\bf Y}_i^@ - {\bf m}(X_{i1}^@)$. So they satisfy the conditions of Theorem 2. Application 
 of Theorem 2 to the two-dimensional stochastic process gives sets $G_{12}$, $G_{11}$, $G_{22}$ due 
 to ordering by time and by intensity and to specific structure of set $A$ in this case.
  
The proof is complete.

\vspace{10 mm}

{\bf Acknowledgement}

The work is supported by Mathematical Center in Akademgorodok 
under agreement No. 075-15-2019-1675 with the Ministry of Science and Higher Education of the Russian Federation.


\bigskip
\begin{thebibliography}{1}



\bibitem{ref01}
Bickel, P. J., and Wichura, M. J., 1971. 
{\em Convergence Criteria for Multiparameter Stochastic Processes and Some Applications}, 
Ann. Math. Stat.{\bf  42} (5), 1656--70. MR0383482


\bibitem{ref0185}
Bouterse, B., and Perros, H., 2012. {\em Scheduling cloud capacity for Time-Varying customer demand}, 2012 IEEE 1st International Conference on Cloud Networking (CLOUDNET), 137-142, doi: 10.1109/CloudNet.2012.6483668.

 
\bibitem{ref025}
Chebunin, M. G., Kovalevskii, A. P., 2021.
{\em  Asymptotics of sums of regression residuals under multiple ordering of regressors}, 
Siberian Electronic Mathematical Reports {\bf 18}, No. 2, 1482--1492.  DOI 10.33048/semi.2021.18.111.
 

\bibitem{ref047}
{Davydov, Y., Egorov, V.}, 2000. 
{\em Functional limit theorems for induced order statistics},
Mathematical Methods of Statistics {\bf 9} (3), 297--313. MR1807096

\bibitem{ref050}
{Davydov, Y., Zitikis, R.}, 2008. 
{\em On weak convergence of random fields},
Annals of the Institute of Statistical Mathematics {\bf 60}, 345--365. MR2403523
%DOI:10.1007/s10463-006-0090-4

\bibitem{ref060}
Dudley, R. M., 1967. 
{\em Measures on non-separable metric spaces},  Illinois Journal of Mathematics
{\bf 11}, 449--453. MR0235087


\bibitem{ref1306}
Ossiander, M., 1987. {\em A Central limit theorem under metric entropy
with $L_2$ bracketing},  Ann. Prob. {\bf 15}, 897--919. MR0893905


\bibitem{ref134}
Straf, M.L., 1972. 
{\em Weak convergence of stochastic processes with several parameters}, Proc. Sixth Berkely Symp.
Math. Statist. Prob. {\bf 2}, 187--222. MR0402847

\end{thebibliography}
\end{document}